\documentclass[11pt]{article}   
\usepackage{amsthm, amsmath}   
\usepackage[english]{babel}

\theoremstyle{plain}
\newtheorem{theorem}{Theorem}[section]
\newtheorem{lemma}[theorem]{Lemma}
\newtheorem{proposition}[theorem]{Proposition}

\theoremstyle{definition}
\newtheorem{definition}[theorem]{Definition}

\newtheorem{algorithm}[theorem]{Algorithm}

\theoremstyle{remark}
\newtheorem{remark}[theorem]{Remark}

\title{A Computational Approach to the ${\mathcal
D}$-module of \\ Meromorphic Functions}

\author{F.J. Castro-Jim\'{e}nez and J.M. Ucha}
\date{September, 2001}

\begin{document}
\maketitle

\begin{abstract} Let $D$ be a divisor in ${\bf C}^n$.
We present methods to compare the ${\mathcal D}$-module of the
meromorphic functions ${\mathcal O}[* D]$ to some natural
approximations. We show how the analytic case can be treated with
computations in the Weyl algebra.
\end{abstract}

\section{Introduction}

Let us denote by ${\cal O} = {\cal O}_{{\bf C}^n}$ the sheaf of
holomorphic functions on $X={\bf C}^n$. Consider a point $p\in
{\bf C}^n$. $Der({\cal O}_p)$ is the ${\cal O}_p$-module of ${{\bf
C}}$-derivations of ${\cal O}_p$. The elements in $Der({\cal
O}_p)$ are called {\it vector fields}.

Let $D\subset X$ be a divisor (i.e. a hypersurface) and $p\in D$.
A vector field $\delta \in Der({\cal O}_p)$ is said to be {\it
logarithmic} with respect to $D$ if $\delta(f)=af$ for some $a\in
{\cal O}_p$, where $f$ is a local (reduced) equation of the germ
$(D,p)\subset ({\bf C}^n,p)$. The ${\cal O}_p$-module of
logarithmic vector fields (or logarithmic derivations) is denoted
by $Der(\log D)_p$. This yields an ${\mathcal O}$-module sheaf
denoted by $Der(\log D)$ (see \cite{Saito}).

Let us denote by ${\cal D} = {\cal D}_X$ the sheaf (of rings) of
linear differential operators with holomorphic coefficients on
$X={\bf C}^n$. A local section $P$ of ${\cal D}$ (i.e. a linear
differential operator) is a finite sum $P = \sum_\alpha a_\alpha
\partial^\alpha$ where $\alpha=(\alpha_1,\ldots,\alpha_n)\in
{\bf N}^n$, $a_\alpha$ is a local section of ${\cal O}$ and
$\partial=(\partial_1,\ldots,\partial_n)$ with
 $\partial_i=\frac{\partial}{\partial x_i}$ in some local chart.

For any divisor $D\subset {\bf C}^n$ we denote by  ${\cal O}[\star
D]$ the sheaf of meromorphic functions with poles along $D$. It
follows from the results of Bernstein-Bj\"ork (\cite{Bernstein},
\cite{Bjork}) on the existence of the $b$-function for each local
equation $f$ of $D$, that ${\mathcal O}[\star D]$ is a left
coherent ${\cal D}$-module. Kashiwara proved that the dimension of
its characteristic variety is $n$ and then that  ${\cal O}[\star
D]$ is holonomic, \cite{Kas-II}.

We will consider some $\mathcal D$-modules associated to any
divisor $D$:
\begin{itemize}
\item The (left) ideal $I^{\log D }\subset \mathcal D$ generated by the
logarithmic vector fields $Der (\log D)$.
\item The (left) ideal $\widetilde{I}^{\log D} \subset {\mathcal D}$
generated by the set $\{ \delta + a\, \vert\, \delta\in I^{\log
D}{\mbox{ and }} \delta(f)=a f \}$. More generally, the ideals
$\widetilde{I}^{(k)\log D}$ generated by the set $$\{ \delta +
ka\, \vert\, \delta\in I^{\log D}{\mbox{ and }} \delta(f)=a f \}$$
\item The modules $M^{\log D}= \mathcal D/I^{\log D}$,
$\widetilde{M}^{\log D} = \mathcal D/\widetilde{I}^{\log D}$ and
more generally $\widetilde{M}^{(k)\log D} = \mathcal
D/\widetilde{I}^{(k) \log D}$.
\end{itemize}

The inclusion $\widetilde{I}^{(k)\log D} \subset Ann_{\mathcal
D}(1/f^k)$ yields to a natural morphism $\phi^k_D :
\widetilde{M}^{(k)\log D} \rightarrow {\mathcal O}[\star D]$
defined by $\phi^k_D(\overline{P}) = P(1/f^k)$ where
$\overline{P}$ denotes the class of the operator $P\in \mathcal D$
modulo $\widetilde{I}^{(k)\log D}$. The image of $\phi^k_D$ is
$\mathcal D\frac{1}{f^k}$, i.e. the ${\mathcal D}$-submodule of
${\mathcal O}[\star D]$ generated by $1/f^k$.

Considering the general ideals $\widetilde{I}^{(k)\log D}$ is a
suggestion of Prof. Tajima. The point is the well known chain of
inclusions $$\mathcal D \cdot {f^{-1}} \subset \mathcal D \cdot
{f^{-2}} \subset \cdots \subset \mathcal D \cdot{f^k} = \mathcal D
\cdot{f^{k-1}} = \cdots = {\mathcal O}[\star D],$$ where $k$ is
least integer root of the $b$-function.

\vspace{.5cm}

We are interested in the germ $(D,p) \subset ({\bf C}^n,p)$ for a
fixed point $p \in D$. So we will work in the ring of germs of
linear differential operators ${\mathcal D}_p$. We suppose that $p
= 0 \in {\bf C}^n$ and from now on, we will denote ${\mathcal D} =
{\mathcal D}_0$. In this context we will use the Weyl algebra
$A_n({\bf C})$ as a subring of ${\mathcal D}$.

Under a computational point of view the divisor $D$ will be
defined by a polynomial $f \in {\bf C}[x_1,\ldots,x_n]$. The
$b$-function of $f$ is computable by \cite{Oaku} and we give a
direct method to present ${\mathcal O}[*D]$. If the calculation is
intractable --as sometimes happens in the examples-- we also
present an indirect method to deduce that ${\mathcal O}[*D]$ and
the modules $\widetilde{M}^{(k)\log D}$ do not coincide. The
method is strongly based in the following result of \cite{MK2}:

\begin{theorem}{\label{ext0}} The restriction to $D$ of the
sheaf $Ext^i_{\mathcal D}({\mathcal O}[*D], {\mathcal O})$ is zero
for $i \geq 0$.
\end{theorem}

More precisely, we prove that, under certain algorithmic
conditions, some cohomology groups are not zero. The interest of
this second method has been tested in \cite{Ucha-Tesis},
\cite{Castro-Ucha-jsc} and \cite{Castro-Ucha-moscu}. We point out
how the algorithms presented in \cite{tsai-walther} --that only
calculate cohomology groups in the algebraic case-- could also be
useful in some {\it analytic} situations.

\vspace{.5cm}

It is important to underline that our methods manage the {\it
analytic} case. As the inclusion $A_n({\bf C}) \subset {\mathcal
D}$ is flat, the computation of syzygies and free resolutions in
the Weyl algebra yields to the analogous computations in
${\mathcal D}$.

\section{Comparison algorithms}{\label{algoritmos}}

We propose in this section two methods to compare the logarithmic
modules presented above. It is important to remark that the
computation of the {\em analytic} $Der(\log D)$ can be made for a
divisor $D$ if its local equation is a polynomial $f \in {\bf
C}[x_1,\ldots,x_n]$ . Simply compute, using Gr\"{o}bner basis, a
system of generators of $(f_1,\ldots,f_n,f)$ where $f_i =
\frac{\partial f}{\partial x_i}$ because the inclusion of the Weyl
algebra in ${\mathcal D}$ is flat.

\subsection{Direct comparison}
The first method is complete but needs the calculation of the
$b$-function.

Experimental evidences show that if the divisor is not locally
{\em Euler homogeneous} (i.e. there is no $\delta \in Der(\log D)$
such that $\delta(f) = f$) the $b$-function is hard to compute.
More precisely, the problem seems to be the calculation of
$Ann_{{\mathcal D}[s]}(f^s)$ and the use of certain elimination
orders during the calculation of Gr\"{o}bner basis (here ${\mathcal
D}[s]$ stands for the polynomial ring, with the indeterminate $s$
commuting with ${\mathcal D}$ ). Nevertheless, the method is
applicable to all the Euler homogeneous divisors we have
considered (and to some not Euler homogeneous).

\begin{algorithm} {\bf INPUT:} A local equation $f = 0$ of a divisor $D
$;
\begin{enumerate}
\item Compute the $b$-function of $f$. Let $-\alpha_0$ be the least
integer root.
\item Compute the ideal $Ann_{\mathcal D}(1/f^{\alpha_0})$.
\item Compute a set of generators $\{ {\bf s}_1, \ldots , {\bf s}_r \}$ of
$Syz(f_1, \ldots,f_n,f)$. The ideal $\widetilde{I}^{(\alpha_0)\log
D}$ is generated by the elements $${\bf s}_j \left(
\begin{array}{c} \partial_1 \\ \vdots \\ \partial_n \\ -\alpha_0
\end{array} \right)\in {\mathcal D}.$$
\item Compare $Ann_{\mathcal D}(1/f^{\alpha_0})$ and  $\widetilde{I}^{(\alpha_0)\log
D}$.
\end{enumerate}
\noindent {\bf OUTPUT:} ${\mathcal O}[*D] \simeq
\widetilde{M}^{(\alpha_0)\log D} \Leftrightarrow Ann_{\mathcal
D}(1/f^{\alpha_0}) = \widetilde{I}^{(\alpha_0)\log D}$.
\end{algorithm}

\vspace{.5cm}

\noindent The correctness of the algorithm is obvious as
$${\mathcal O}[*D] \simeq {\mathcal D}\frac{1}{f^{\alpha_0}}
\simeq {\mathcal D}/Ann_{\mathcal D}(1/f^{\alpha_0}).$$

\subsection{Indirect deduction: a sufficient condition}

This second method is an alternative way when you can not obtain
the $b$-function. In the worst case, it only needs the computation
of a free resolution of $\widetilde{M}^{(\alpha)\log D}$ for any
integer $\alpha \geq 1$. More precisely, the algorithm looks for a
technical condition in some step of the free resolution. In many
examples, it is enough to compute only the first
syzygies\footnote{Taking into account that computing a complete
free resolution can be a problem of great complexity, this option
is very interesting.}.

\begin{definition}If $$ 0 \rightarrow {\mathcal D}^{r_s} \stackrel{\varphi_s}
\longrightarrow \cdots \rightarrow {\mathcal D}^{r_2}
\stackrel{\varphi_2}\longrightarrow {\mathcal D}^{r_1}
\stackrel{\varphi_1}\longrightarrow {\mathcal D}
\stackrel{\pi}\longrightarrow M \rightarrow 0$$ is a free
resolution of a module $M$, we will say that the {\em Successive
Matrices Condition (SMC) holds at level $i$} if the two succesive
morphisms $\varphi_{i},\varphi_{i+1}$ have matrices verifying:
\begin{enumerate} \item The matrix of $\varphi_{i+1}$ has
no part in ${\mathcal O}$ in some column $j$, i.e. all the
elements in the $j$-th column are in the (left) ideal generated by
$\partial_1, \ldots, \partial_n$.
\item The matrix of $\varphi_{i}$ has no constants in the in the $j$-th
row. That is, for each $P$ in the $j$-th row $P(1)$ is a function
$h$ with $h(0) = 0$.
\end{enumerate}
\end{definition}

\begin{algorithm}{\label{indirect}} {\bf INPUT:} A local equation $f = 0$ of a divisor $D
$;
\begin{enumerate}
\item Compute a set of generators $\{ {\bf s}_1, \ldots , {\bf s}_r \}$ of
$Syz(f_1, \ldots,f_n,f)$. The ideal $\widetilde{I}^{(\alpha)\log
D}$ is generated by the elements $${\bf s}_j \left(
\begin{array}{c} \partial_1 \\ \vdots \\ \partial_n \\ -\alpha
\end{array} \right)\in {\mathcal D}.$$
\item Compute a free resolution of $M = \widetilde{M}^{(\alpha)\log
D}$ $$ 0 \rightarrow {\mathcal D}^{r_s} \stackrel{\varphi_s}
\longrightarrow \cdots \rightarrow {\mathcal D}^{r_2}
\stackrel{\varphi_2}\longrightarrow {\mathcal D}^{r_1}
\stackrel{\varphi_1}\longrightarrow {\mathcal D}
\stackrel{\pi}\longrightarrow M \rightarrow 0.$$
\end{enumerate}
\noindent {\bf OUTPUT:} {\bf IF} SMC holds {\bf THEN} $ {\mathcal
O}[*D] \neq \widetilde{M}^{(\alpha)\log D}.$
\end{algorithm}

\vspace{.5cm} We need a lemma to justify the algorithm. It
explains the role of the SMC. The idea is obtaining an element in
$Ker \varphi_{i+1}$ that is not in $Im \varphi_{i}$.

\begin{lemma}{\label{extnocero}} Let $D$ be a divisor and
$$ 0 \rightarrow {\mathcal D}^{r_s} \stackrel{\varphi_s}
\longrightarrow \cdots \rightarrow {\mathcal D}^{r_2}
\stackrel{\varphi_2}\longrightarrow {\mathcal D}^{r_1}
\stackrel{\varphi_1}\longrightarrow {\mathcal D}
\stackrel{\pi}\longrightarrow \widetilde{M}^{(\alpha)\log D}
\rightarrow 0 \ \  (*) $$ a free resolution of
$\widetilde{M}^{(\alpha)\log D}$ that verifies SMC at level $i$.
Then $$Ext^{i}_{\mathcal D}(\widetilde{M}^{(\alpha)\log
D},{\mathcal O}) \neq 0.$$
\end{lemma}
\begin{proof}
To obtain the $Ext$ groups, we have to apply the functor
$Hom_{\mathcal D}(-, {\mathcal O})$ to $(*)$. Using that
$$Hom_{\mathcal D}({\mathcal D}^r ,{\mathcal O}) \simeq {\mathcal
O}^r$$ we obtain the complex $$ 0 \rightarrow {\mathcal O}
\stackrel{\varphi_1^t}\longrightarrow {\mathcal O}^{r_1}
\stackrel{\varphi_2^t}\longrightarrow {\mathcal O}^{r_2}
\rightarrow \cdots \stackrel{\varphi_s^t}\longrightarrow {\mathcal
O} \stackrel{\varphi_s^t}\longrightarrow {\mathcal O}^{r_s}
\rightarrow 0,$$ where $\varphi_i^t$ denotes the morphism with
matrix the transposed of $\varphi_i$. The derivatives now act
naturally.

Then $$Ext^{i}_{\mathcal D}(\widetilde{M}^{(\alpha)\log
D},{\mathcal O}) = Ker \varphi_{i+1}^t/ Im \varphi_{i}^t.$$ If the
matrix of $\varphi_{i+1}$ has no part in ${\mathcal O}$ in the
$j$-th row, then ${\bf e} = (0,\ldots,1,\ldots,0)$ --where $1$ is
in the $j$-th position-- is in $Ker \varphi_{i+1}^t$, as the
derivatives applied to 1 are zero.

This element can not be in $Im \varphi_{i}^t$ if the matrix has no
constants. Applying the operators of the matrix it is not possible
to obtain elements of degree 0.
\end{proof}

We have the key to state the main result of this section: the
correctness of \ref{indirect}.

\begin{proposition} Let $D$ be a divisor with a free resolution of
$\widetilde{M}^{(\alpha)\log D}$ that verifies SCM at some level.
Then $ {\mathcal O}[*D] \neq \widetilde{M}^{(\alpha)\log D}.$
\end{proposition}
\begin{proof} Evident from \ref{extnocero} and \ref{ext0}.
\end{proof}

\vspace{.5cm}

A special case of SCM appears when you have a resolution of type
$$ 0 \rightarrow {\mathcal D}^{r_n} \stackrel{\varphi_n}
\longrightarrow \cdots \rightarrow {\mathcal D}^{r_2}
\stackrel{\varphi_2}\longrightarrow {\mathcal D}^{r_1}
\stackrel{\varphi_1}\longrightarrow {\mathcal D}
\stackrel{\pi}\longrightarrow M \rightarrow 0 \ $$ of length $n$.
Then $$Ext^n_{\mathcal D}(M,{\mathcal O}) = \frac{{\mathcal
O}^{r_n}}{Im \varphi_n^t},$$ and SCM means, at level $n$, that you
can find in the matrix of $\varphi_n$ a row with no constants.

\vspace{.5cm}

\begin{remark}
Of course, natural generalizations of the SCM condition has to do
with finding explicit elements in some $Ker \varphi_{i+1}^t$ with
special properties. It is not that easy in general! Nevertheless
the results of \cite{tsai-walther} can be applied in this
situation as follows: $$Ext_{A_n({\bf C})}^0(R[*D],R) \neq 0
\Rightarrow Ext_{{\mathcal D}}^0({\mathcal O}[*D],{\mathcal O})
\neq 0,$$ where $R$ is the ring of polynomials and $R[*D]$ is its
localization with respect to the equation $f$ of $D$.
\end{remark}

\section{Application to the Spencer case.}{\label{aplicaciones}}

In this section we explain how to apply the sufficient condition
to a special case in which a tailored free resolution is provided.

\vspace{.5cm}

\begin{definition}\label{free} {\rm (\cite{Saito})} The divisor $D$ is said to
be {\it free at the point} $p\in D$ if the ${\cal O}_{p}$-module
$Der(\log D)_p$ is free. The divisor $D$ is called {\it free} if
it is free at each point $p\in D$.
\end{definition}

Smooth divisors and normal crossing divisors are free. By
\cite{Saito} any reduced germ of plane curve $D\subset {\bf C}^2$
is a free divisor.

By Saito's criterium \cite{Saito}, $D\equiv(f=0)\subset {\bf C}^n$
is free at a point $p$ if and only if there exist $n$ vector
fields $\delta_i = \sum_{j=1}^n a_{ij}
\partial_j$, $i=1,\ldots,n$, such that $\det(a_{ij})=uf$ where
$u$ is a unit in ${\mathcal O}_p$. Here $\partial_j$ is the
partial derivative $\frac{\partial}{\partial x_j}$ and $a_{ij}$ is
a holomorphic function in ${\mathcal O}_p$.

\begin{definition}
We say that a free divisor $D$ is {\it of Spencer type} if the complex $$
{\mathcal D}\otimes_{\mathcal O} \wedge^\bullet Der(\log D) \rightarrow M^{\log
D} \rightarrow 0$$ (introduced in \cite{Cald3}) is a (locally) free resolution
of $M^{\log D}$ and if this last ${\mathcal D}$-module is holonomic.
\end{definition}

There are analogous resolutions for the family of modules
$\widetilde{M}^{(k)\log D}.$

\vspace{.5cm}

For this family of divisors, the solution complex $Sol(M^{\log
D})$ (that is, the complex ${\bf R}{\mathcal H}om_{\mathcal
D}(M^{\log D},{\mathcal O})$) is naturally quasi-isomorphic to
$\Omega^\bullet(\log D)$ (as we pointed in
\cite{Castro-Ucha-moscu} as a deduction of \cite{Cald3}). On the
other hand, a duality theorem proved in \cite{Castro-Ucha-moscu}
has important consequences comparing $\widetilde{M}^{\log D}$ and
${\mathcal O}[\star D]$, namely

\begin{theorem} {\rm (\cite{Ucha-Tesis,Castro-Ucha-jsc})} In dimension 2,
the morphism $\phi^1_D$ is an isomorphism if and only if $D$ is a
quasi-homogeneous plane curve.
\end{theorem}

\begin{theorem} {\rm \cite{Castro-Ucha-moscu}} Suppose the
divisor $D\subset {\bf C}^n$ is free and locally
quasi-homogeneous. Then the morphism $\phi^1_D$ is an isomorphism
(so, $\widetilde{M}^{\log D}$ and ${\mathcal O}[\star D]$ are
isomorphic as ${\mathcal D}$-modules).
\end{theorem}

\vspace{0.5cm}

The methods presented in section \ref{algoritmos} give us
computational tools to check the comparison between
$\widetilde{M}^{\log D}$ and ${\mathcal O}[\star D]$ .

\begin{remark} Once you have the duality of \cite{Castro-Ucha-moscu}, you also have
a strategy to study the {\em Logarithmic Comparison Theorem
(LCT)}, that is, the complex $\Omega^\bullet(\star D)$ of
meromorphic differential forms and the complex
$\Omega^\bullet(\log D)$ are quasi isomorphic. You have to travel
round the following chain of isomorphisms: $$\Omega^\bullet(\star
D) \simeq DR({\mathcal O}[\star D]) \simeq DR({\mathcal
D}/Ann_{\mathcal D}(1/f)) \simeq DR(\widetilde{M}^{\log D})
\simeq$$ $$\simeq DR((M^{\log D})^*) \simeq Sol(M^{\log D}) \simeq
\Omega^\bullet(\log D),$$ where for each coherent ${\mathcal
D}$--module $M$ we denote by $DR(M)$ the de Rham complex of $M$
(see \cite{MK2}).
\end{remark}

\begin{remark} There are two interesting experimental suggestions:
\begin{itemize}
\item We don't know examples of free divisors with integer roots of
their $b$-function less than -1. \item We only know free divisors
of Spencer type.
\end{itemize}
\end{remark}

\vspace{.5cm}

Finally, we have the following result:

\begin{proposition}
Let $D\equiv(f=0)$ be a Spencer divisor. Let
${\delta_1,\ldots,\delta_n}$ be a basis of $Der(\log D)$ with
$\delta_i = \sum_{j=1}^n{a_ij}\partial_j$ for $1 \leq j \leq n$.
If $\sum_{j=1}^n\partial_j (a_{ij})= 0$, then $\widetilde{M}^{\log
D}$ and ${\mathcal O}[\star D]$ are not isomorphic.
\end{proposition}
\begin{proof}
The last matrix of the Spencer free resolution of
$\widetilde{M}^{\log D}$ is of a very special type . Its elements
(due to the duality formulas of \cite{Castro-Ucha-moscu}) are of
the form $$\delta_i + \sum_{j=1}^n{\partial_j (a_{ij})}.$$ So,
applying lemma \ref{extnocero} $Ext^n_{\mathcal
D}(\widetilde{M}^{\log D}, {\mathcal O}) \neq 0,$ so
$\widetilde{M}^{\log D}$ and ${\mathcal O}[\star D]$ are not
isomorphic. Thus there is no possible quasi-isomorphism between
$DR({\mathcal D}/Ann_{\mathcal D}(1/f))$ and
$DR(\widetilde{M}^{\log D})$.
\end{proof}

\section{Examples}

In the following examples, the computation of syzygies among
polynomials have been made with {\tt CoCoA} (see \cite{CoCoA}).

The computations of syzygies in the Weyl Algebra, global
$b$-functions and ideals of type $Ann_{\mathcal D}(1/f^{\alpha})$
have been made with {\tt kan/sm1}, \cite{Kan}, that is, using the
algorithms of \cite{Oaku}.

\subsection{Example 1: $D\equiv(x (x^2 - y^3)(x^2 - z y^3) = 0)$}

We will treat here the divisor $D \subset {\bf C}^3$ whose local
equation at $(0,0,0)$ is given by $f = 0$ with $$f = x (x^2 - y^3)
(x^2 - z y^3).$$

This divisor is (globally) free and $\delta_1, \delta_2, \delta_3$
form a (global) basis of $Der (log D)$, where $$\begin{array}{ccl}
\delta_1 & = & \frac{3}{2} x
\partial_x + y \partial_y \\ \delta_2 & = & (y^3z - x^2)
\partial_z \\ \delta_3 & = & (-\frac{1}{2}xy^2) \partial_x - \frac{1}{3} x^2
\partial_y + (y^2z^2- y^2z) \partial_z,
\end{array}$$
\noindent whose coefficients verify that $$\left|
\begin{array}{ccl} \frac{3}{2} x & y & 0 \\ 0 & 0 &
y^3z - x^2 \\ -\frac{1}{2}xy^2 & - \frac{1}{3} x^2 & y^2z^2- y^2z
\end{array}
\right| = -\frac{1}{2} f.$$

\vspace{.5cm}

To begin with, we have to follow two steps:
\begin{itemize}
\item Step 1: Verify that $M^{\log D}$ is holonomic\footnote{This computation
could be made with \cite{Kan}}. The interest of this question is
evident: if $M^{\log D}$ is not holonomic, the computation of its
dual could not be managed as we do.

\item Step 2: Compute a free resolution of $M^{\log D}$
with Gr\"{o}bner basis computation of syzygies.  Check if $D$ has a
free resolution of Spencer type. If this happens then duality
holds by \cite{Castro-Ucha-moscu}.
\end{itemize}

The example verifies these properties:

\begin{enumerate}
\item The module
$Syz (\delta_1, \delta_2, \delta_3)$ is generated by the syzygies
obtained from the commutators $[\delta_i, \delta_j]$. We have $Syz
(\delta_1, \delta_2 \delta_3) = \langle {\bf s}_{12}, {\bf s}_{13}
, {\bf s}_{23} \rangle$ where $$
\begin{array}{ccl}
{\bf s}_{12} & = & (-\delta_2, \delta_1 - 3, 0) \\ {\bf s}_{13} &
= & (-\delta_3, 0, \delta_1 - 2) \\ {\bf s}_{12} & = & (0,
-\delta_3 - y^2z, \delta_2).
\end{array}
$$
\item On the other hand, the module $Syz ({\bf s}_{12}, {\bf s}_{13},
{\bf s}_{23})$ is generated by the element ${\bf r}$: $$ {\bf r} =
(-y^2z^2\partial_z + y^2z\partial_z + \frac{1}{2} xy^2 \partial_x
- y^2z + \frac{1}{3} x^2
\partial_y,\ \ y^3z
\partial_z - x^2 \partial_z, \ \ -y \partial_y - \frac{3}{2} x
\partial_x + 5).$$
This is the element required to have the Spencer type resolution
so, as we have said, duality holds.
\end{enumerate}

\noindent We calculate the $b$-function of $f$. Its least integer
root is -1, so $${\mathcal O}[\star D] \simeq {\mathcal D}\cdot
1/f \simeq Ann_{\mathcal D}(1/f).$$ To finish, we check that
$\widetilde{I}^{\log D} = Ann_{\mathcal D}(1/f)$.

\subsection{Example 2: $D\equiv(x^4 + y^4 + z^4 + x^2y^2z^2 = 0)
\subset {\bf C}^3$}{\label{nolibre}}

This is divisor is not free. A set of generators of
$\widetilde{M}^{\log D}$ is $$\begin{array}{ccl} \delta_1 & = &
(1/8x^2y^3z^2 + y) + ( -1/32x^3y^3z^2 - 1/16xyz^4 - 1/4xy)
\partial_x +
\\
 & & + (-1/32x^2y^4z^2 + 1/8x^2z^2 - 1/4y^2) \partial_y + (-1/4yz )
\partial_z,\\
 & & \\
\delta_2 & = & (-1/16x^4y^4z + z) +  (1/64x^5y^4z + 1/32x^3y^2z^3
- 1/4xz) \partial_x +
\\ & & + (1/64x^4y^5z - 1/16x^4yz - 1/4yz) \partial_y + (1/8x^2y^2 -
 1/4z^2)\partial_z,\\
  & & \\
\delta_3 &  = &  (-1/16x^4y^3z^3 - 1/8x^2y^5z) + \\ & & +
(1/64x^5y^3z^3 + 1/32x^3y^5z + 1/32x^3yz^5 + 1/16xy^3z^3)
\partial_x + \\
 & & + (1/64x^4y^4z^3 + 1/32x^2y^6z - 1/16x^4z^3 - 1/8x^2y^2z - 1/4z^3) \partial_y
 + \\ & & +
 (1/8x^2yz^2 + 1/4y^3) \partial_z,\\
 & & \\
\delta_4 & = & (-1/16xy^4z^4 + x) + (1/64x^2y^4z^4 + 1/32y^2z^6 +
1/8y^2z^2 - 1/4x^2)
\partial_x + \\ & &  + (1/64xy^5z^4 - 1/16xyz^4 - 1/4xy) \partial_y +(-1/4xz)
\partial_z,\\
 & & \\
\delta_5 & = & (1/8x^2y^2z + 1/4z^3) \partial_x +(-1/8xy^2z^2 -
1/4x^3) \partial_z, \\
 & & \\
\delta_6 & = & (1/8xy^5z^2 + 1/4x^3y^3) + \\ & & + (-1/32x^2y^5z^2
- 1/16x^4y^3 - 1/16y^3z^4 - 1/8x^2yz^2 - 1/4y^3) \partial_x + \\
 & & + (-1/32xy^6z^2 - 1/16x^3y^4 + 1/8xy^2z^2 + 1/4x^3)
\partial_y. \end{array}$$

The free resolution is huge, but anyway computable with {\tt
kan/sm1}. It is of type $$ 0 \rightarrow {\mathcal D}
\stackrel{\varphi_3}\longrightarrow {\mathcal D}^{r}
\stackrel{\varphi_2}\longrightarrow {\mathcal D}^{s}
\stackrel{\varphi_1}\longrightarrow {\mathcal D}
\stackrel{\pi}\longrightarrow M \rightarrow 0 \ $$ To use
proposition \ref{extnocero} you can check that there the elements
in the last matrix $\varphi_3$ has no constants. So
$$Ext^3_{\mathcal D}(\widetilde{M}^{\log D}, {\mathcal O}) =
{\mathcal O}/ im(\varphi_3)^t \neq 0.$$ You have ${\mathcal O}[*
D] \neq \widetilde{M}^{\log D}$.

\subsection{Example 3: $D\equiv(x^4+y^5+xy^4+zx^6 = 0) \subset
{\bf C}^3$}

In this example we detect that we have, in fact, a (non-Euler
homogeneous) product that is a free divisor\footnote{By the way,
this situation is impossible in dimension 2.}. The basis of
$Der(\log D)$ is $$\begin{array}{ccl} \delta_1 & = &
(x^6z^2 + 5/4x^5yz^2 - 2x^4z - x^2 - 5/4xy)\partial_x + \\
 & & + (5/4x^5yz^2 + 3/2x^4y^2z^2 - 5/2x^3yz - 1/2x^2y^2z - 3/4xy -
 y^2) \partial_y + \\
 & & + (2x^3z^2 - 2xz) \partial_z, \\
  & & \\
\delta_2 & = &
(-4/5x^4yz - 8/5x^3y^2z - 3/4x^2y^3z + 25/4x^4z + 125/16x^3yz -
xy^2 - 1/4y^3 + \\ & & +  125/16xy) \partial_x + (-x^3y^2z -
39/20x^2y^3z - 9/10xy^4z + 125/16x^3yz + 75/8x^2y^2z - \\ & & -
3/4y^3 + 1/4x^2 - 5/16xy + 25/4y^2) \partial_y + (-8/5xyz -
6/5y^2z + 25/2xz) \partial_z, \\
 & & \\
\delta_3 & = &
(-3/10x^3y^2z -3/8x^2y^3z + 25/8x^4z + \\ & & + 125/32x^3yz +
8/25x^3 - 1/8y^3 + 125/32xy) \partial_x + (-3/8x^2y^3z - 9/20xy^4z
+ \\
 & & + 125/32x^3yz + 75/16x^2y^2z + 2/5x^2y -
1/50xy^2 + 1/40y^3 + 1/8x^2 - 5/32xy + \\ & & + 25/8y^2)
\partial_y + (-3/5y^2z + 25/4xz + 16/25) \partial_z. \end{array}$$

>From the monomial $16/25 \partial_z$ we deduce by the Flow Theorem
that $D$ is a product. It is of Spencer type and the third matrix
can be used to show that $Ext^3_{\mathcal D}(\widetilde{M}^{\log
D},{\mathcal O}) \neq 0$ so ${\mathcal O}[* D] \neq
\widetilde{M}^{\log D}$.

\subsection{Example 4: $D\equiv((x+y)(xz+y)(x^4 + y^5 + xy^4) = 0)
\subset {\bf C}^3$}

In this last example the divisor $D \subset {\bf C}^3$ has as a
local equation at $(0,0,0)$ the form $$f = (x+y)(xz+y)(x^4 + y^5 +
xy^4) = 0.$$ \noindent The divisor is globally free with
$\delta_1, \delta_2, \delta_3$ as a global basis of $Der (log D)$.
$$\begin{array}{ccl} \delta_1 & = & (-x^2 - 5/4xy)
\partial_x +(-3/4xy - y^2) \partial_y + (-1/4xz^2 + 1/4xz)
\partial_z \\ & & \\ \delta_2 & = & (xz + y) \partial_z \\ & & \\ \delta_3 & = &
(2x^2y^2 + 5/2xy^3 + 1/2y^4 + 5/2x^3 - 7x^2y - 35/4xy^2 - 11x^2 -
55/4xy) \partial_x + \\ & & + (3/2xy^3 + 3/2y^4 - 1/2x^3 + 2x^2y -
21/4xy^2 - 7y^3 - 33/4xy - 11y^2) \partial_y +
\\ & & + (1/2xy^2z^2 + 1/2y^3z^2 - 1/2xy^2z - 1/2y^3z - 1/4xyz^2 +
 7/4xyz + 3/2y^2z + \\ & & + 3/4xz^2 + 1/2x^2 + 1/2xy + 5/2xz + 7/2yz -
 1/4y) \partial_z.
\end{array}$$

In this case, $D$ is of Spencer type and the third matrix has no
constants so, again we have that $Ext^3_{\mathcal
D}(\widetilde{M}^{log D},{\mathcal O})\neq 0$.

\section{Acknowledgements}
We thank Prof. Tajima and Prof. David Mond for very helpful ideas
and comments.

\end{document}